\documentclass[11pt]{article}

\usepackage{amsmath,amssymb,amsthm,mathtools}

\numberwithin{equation}{section}

\newtheorem{theorem}{Theorem}[section]
\newtheorem{proposition}[theorem]{Proposition}

\newtheorem{remark}[theorem]{Remark}

\newcommand{\R}{\mathbb{R}}

\begin{document}

\title{Large data global well-posedness for a one-dimensional quasilinear wave equation}
\author{
Yuusuke Sugiyama\\
Department of Mathematics, Tokyo University of Science\\
Kagurazaka 1-3, Shinjuku-ku, Tokyo 162-8601, Japan\\
\texttt{sugiyama.y@rs.tus.ac.jp}
}
\date{}
\maketitle

\begin{abstract}
In this paper, we prove global well-posedness with large initial data for the one-dimensional
quasilinear wave equation
\begin{equation*}
  u_{tt}=c(u)^2u_{xx},
  \qquad (t,x)\in (0,T)\times\R,
\end{equation*}
where \(c\) is a positive, bounded, monotonically increasing function with
bounded derivative. This result  gives a
partial resolution of an open problem posed by Glassey, Hunter and Zheng on
the global existence of smooth solutions to this equation for large initial data. Our proof is based on upper
and lower estimates for the Riemann variables via a new comparison principle.
\end{abstract}

\section{Introduction}

\subsection{Main theorem}

In this paper, we consider the Cauchy problem for the one-dimensional
quasilinear wave equation
\begin{equation}\label{eq:main}
  u_{tt}=c(u)^2u_{xx},
  \qquad (t,x)\in (0,T)\times\R,
\end{equation}
with initial data
\begin{equation}\label{eq:initial}
  u(0,x)=u_0(x),\qquad
  u_t(0,x)=u_1(x),
  \qquad x\in\R .
\end{equation}
Here \(u=u(t,x)\) is a real-valued unknown function, and
\(c=c(\theta)\) is a given smooth function.

Throughout this paper, we assume that for  $\theta\in\R$,
\begin{equation}\label{ass:c}
  c\in C^\infty(\R),\qquad
  0<c_*\leq c(\theta)\leq c^*<\infty,
  \qquad
  c'(\theta)\geq 0
\end{equation}
and
\begin{equation}\label{ass:cp}
  \|c'\|_{L^\infty(\R)} = \sup_{\theta \in \R}c'(\theta) <\infty.
\end{equation}

The main result of this paper is the following.

\begin{theorem}\label{thm:main}
Let \(s>1/2\). Assume \eqref{ass:c} and \eqref{ass:cp}. Then, for any
\[
  (u_0,u_1)\in H^{s+1}(\R)\times H^s(\R),
\]
the Cauchy problem \eqref{eq:main}--\eqref{eq:initial} has a unique global
solution
\[
  u\in \bigcap_{j=0}^2 C^j([0,\infty);H^{s+1-j}(\R)).
\]
\end{theorem}

The following functions satisfy the assumptions on \(c\).

\[
c(\theta)=c_0+\delta \tanh \theta \quad (c_0>\delta>0).
\]

\[
c(\theta)=c_-+(c_+-c_-)\frac{1}{1+e^{-\theta}} \quad (0<c_-<c_+).
\]

\[
c(\theta)=c_0+\delta \arctan \theta \quad \left(c_0>\frac{\pi}{2}\delta>0\right).
\]

More generally
\[
c(\theta)=c_-+(c_+-c_-)H(\theta),
\]
where \(H\in C^\infty(\mathbb{R})\), \(0\leq H\leq1\), \(H'\geq0\), and \(\|H'\|_{L^\infty(\mathbb{R})}<\infty\).

\subsection{Known results}
The equation \eqref{eq:main} is the case \(\lambda=0\) of the following
parameterized family of one-dimensional quasilinear wave equations:
\begin{equation}\label{eq:parameterized}
  u_{tt}
  =
  c(u)^2u_{xx}
  +
  \lambda c(u)c'(u)(u_x)^2,
  \qquad 0\leq \lambda\leq 2.
\end{equation}
This family was introduced and studied by Glassey, Hunter and Zheng
\cite{GHZ2}, and was also discussed by Chen and Shen \cite{GCYS}.
Glassey, Hunter and Zheng in \cite{GHZ2} pointed out that, since this equation
does not contain the lower-order term proportional to \(u_x^2\), it appears
possible that it has global smooth solutions for arbitrary smooth initial
data, but that a complete proof or a counterexample was open.
The present paper gives an affirmative answer to this problem under the
additional assumptions that \(c\) is uniformly positive, bounded,
monotonically increasing, and has bounded derivative.

If \(\lambda=2\), then the parameterized nonlinear wave equation
\eqref{eq:parameterized} is formally equivalent to the conservation system
\begin{equation*}
  \partial_t
  \begin{pmatrix}
    U \\ V
  \end{pmatrix}
  -
  \partial_x
  \begin{pmatrix}
    V \\ p(U)
  \end{pmatrix}
  =0,
\end{equation*}
where
\[
  U(t,x)=u(t,x),\qquad
  V(t,x)=\int_{-\infty}^{x} u_t(t,y)\,dy,
  \qquad
  p'(\theta)=c(\theta)^2 .
\]
This conservation system is referred to as a \(p\)-system and describes
several phenomena of wave propagation in nonlinear media, including
electromagnetic waves in transmission lines, shearing motion in
elastic-plastic rods, and one-dimensional gas dynamics.
In the study of the \(p\)-system, the following variables, which are related
to the derivatives of the Riemann invariants, play an important role:
\begin{equation} \label{Riev}
  R=u_t+c(u)u_x,\qquad
  S=u_t-c(u)u_x .
\end{equation}
For hyperbolic systems of conservation laws, necessary and sufficient
conditions for the global existence of classical solutions have been
extensively studied. Roughly speaking,  if
\(R\) and \(S\) are non-positive at the initial time, then the corresponding
classical solution exists globally in time. On the other hand, if either
\(R\) or \(S\) takes a positive value at some point initially, then a
singularity, such as a shock, may form in finite time (e.g. Johnson \cite{JJ}, Klainerman and  Majda \cite{km}, Lax \cite{lax} and  Chen, Pan and Zhu \cite{GPZ} ).

When \(\lambda=1\), the equation in \eqref{eq:parameterized} is called the
variational wave equation:
\[
  u_{tt}-c(u)\bigl(c(u)u_x\bigr)_x=0.
\]
As its name suggests, this equation has a variational structure. More precisely,
it is derived from the least action principle
\[
  \frac{\delta}{\delta u}
  \int \left\{ u_t^2-c(u)^2u_x^2 \right\}\,dx\,dt=0.
\]
Moreover, it possesses the energy conservation law
\[
  E(t)
  =
  \int_{\R}
  \left\{
    u_t(t,x)^2+c(u(t,x))^2u_x(t,x)^2
  \right\}\,dx
  =
  E(0).
\]
This conservation law can be regarded as a quasilinear analogue of the standard
energy conservation law for the linear wave equation.
The variational wave equation has physical backgrounds including nematic
liquid crystals and long waves on a dipole chain in the continuum limit
(see \cite{GHZ2}).
Zhang and Zheng \cite{zz1} showed the existence of global classical solutions
under the assumption that
\[
  R(0,x)\leq0,\qquad S(0,x)\leq0
  \quad (x\in\R).
\]
Glassey, Hunter and Zheng \cite{GHZ,GHZ2} constructed finite-time blow-up
solutions for suitable initial data violating this non-positive condition.

Using the method of Zhang and Zheng \cite{zz1}, the author in \cite{s1}
extended their global existence result to the range \(0\leq\lambda\leq2\)
under the same non-positive condition on the initial Riemann variables.
Under this condition, the non-positivity of \(R\) and \(S\) is preserved;
that is,
\[
  R(t,x)\leq0,\qquad S(t,x)\leq0
  \qquad (t\geq0,\ x\in\R).
\]
Moreover, one can obtain time-independent lower bounds for \(R\) and \(S\).
In particular, there exists a constant \(C>0\), depending only on the initial
data and \(c\), such that
\[
  -C\leq R(t,x)\leq0,\qquad
  -C\leq S(t,x)\leq0
  \qquad (t\geq0,\ x\in\R).
\]
In \cite{s6}, the blow-up result for the variational case \(\lambda=1\) was
extended to the range \(0<\lambda\leq1\).
The construction of blow-up solutions for the case \(1<\lambda<2\) had
remained open.

In \cite{s3,s4}, the author studied the degeneracy of solutions to
\eqref{eq:parameterized} for \(0\leq \lambda\leq2\), in the case where
\(c'(\theta)\geq0\) but \(c\) is not assumed to be uniformly positive from
below. Typical examples include \(c(u)=1+u\). In this setting, even if the
initial wave speed satisfies
\[
  c(u_0(x))\geq \delta>0,
  \qquad x\in\R,
\]
it may happen that
\[
  \inf_{x\in\R} c(u(t,x))\to0
\]
in finite time. This means that the strict hyperbolicity of the equation
breaks down. The papers \cite{s3,s4} investigated sufficient conditions for
the occurrence of such degeneracy.
Using the variables introduced above in \eqref{Riev}, the equation
\eqref{eq:parameterized} can be rewritten as the first-order system
\begin{equation}\label{fsin}
 \left\{
\begin{aligned}
\partial_t R -c(u)\partial_x R
&=
\dfrac{c'(u)}{4c(u)}
\left\{
 \lambda R^2 + 2(1-\lambda)RS -(2-\lambda)S^2
\right\}, \\
\partial_t S +c(u)\partial_x S
&=
\dfrac{c'(u)}{4c(u)}
\left\{
 \lambda S^2 + 2(1-\lambda)RS -(2-\lambda)R^2
\right\}.
\end{aligned}
\right.
\end{equation}
When \(\lambda=0\), the potentially explosive quadratic terms
\(\lambda R^2\) and \(\lambda S^2\) disappear from the right-hand side.
For this reason, it has been expected that solutions should not develop
gradient blow-up and that global well-posedness should hold in this case.

On the other hand, in the construction of blow-up solutions for
\(0<\lambda\leq1\), one of the main difficulties is to control the
dissipative terms
\[
  -(2-\lambda)S^2
  \qquad\text{and}\qquad
  -(2-\lambda)R^2,
\]
which may prevent the Riccati-type growth caused by
\(\lambda R^2\) or \(\lambda S^2\).
Thus, the treatment of these negative quadratic terms is a key point in the
construction of finite-time blow-up solutions.

\subsection{Strategy of the proof and the idea of the comparison principle}

We first prove the theorem for smooth initial data
\[
  (u_0,u_1)\in H^3(\mathbb{R})\times H^2(\mathbb{R}).
\]
In this case, the corresponding local solution is sufficiently regular and
satisfies the equation in the classical sense. Hence we may use the Riemann variables introduced in \eqref{Riev} and rewrite
the equation as a first-order system for \(R\) and \(S\).

The proof is based on pointwise estimates for the Riemann variables along
the characteristic curves. The monotonicity on the characteristic curves of the weighted quantities
\[
  \frac{R}{\sqrt{c(u)}},
  \qquad
  \frac{S}{\sqrt{c(u)}}
\]
gives the uniform upper bound
\[
  R(t,x)\leq P,\qquad S(t,x)\leq P,
\]
where \(P\) depends only on the initial data and the uniform upper and lower
bounds of \(c\).
The main point is to obtain lower estimates for \(R\) and \(S\). For this
purpose, we compare them with the solution of the ordinary differential
equation
\[
  y'(t)=A\{Py(t)-P^2\},
  \qquad
  y(0)=m_0,
\]
where \(m_0\) is determined by the lower bounds of the initial Riemann variables
\(R(0,\cdot)\) and \(S(0,\cdot)\)
and
\[
A=\left\|\frac{c'}{2c}\right\|_{L^\infty} =\sup_{\theta \in  \R} \frac{c'(\theta)}{2c (\theta)}.
\]
The comparison
principle shows that
\[
  R(t,x)\geq y(t),
  \qquad
  S(t,x)\geq y(t).
\]
Combining the upper and lower estimates, we obtain
\[
  \sup_{0\leq t\leq T}
  \left(
    \|R(t)\|_{L^\infty}
    +
    \|S(t)\|_{L^\infty}
  \right)
  <\infty
\]
for every finite \(T<T^*\).
Here we explain how the ordinary differential equation for
\(y\) naturally arises from the equations for the Riemann variables.
We consider 
\[
  m(t)=\min\left\{
    0,\inf_{x\in\mathbb{R}}R(t,x),
      \inf_{x\in\mathbb{R}}S(t,x)
  \right\}.
\]
Formally, we want to show that
\[
  m'(t)\geq A\{P m(t)-P^2\}.
\]
For the moment, we suppose that
the infimum is attained and that
\[
  m(t)=R(t,x_t) \leq 0.
\]
At this point, we have
\[
  R(t,x_t)=m(t),
  \qquad
  m(t)\leq S(t,x_t)\leq P.
\]
By considering the first equation in \eqref{fsin} with $\lambda =0$ at $x=x_t$, since $R_x(t,x_t)=0$, we formally have 
\[
  m'(t)
  =
  R_t (t,x_t) - c(u) R_x (t,x_t)
  =
  \frac{c'(u)}{2c(u)}
  S(t,x_t)\bigl(m(t)-S(t,x_t)\bigr).
\]
For $S(t,x_t) \in [m(t), P]$, we have
\[
m'(t) \geq A(m(t) -P)P,
\]
which derives the differential equation for $y$.

Upper and lower estimates for $R$ and $S$ yield that
the blow-up criterion for the local solution then excludes finite-time
breakdown. This proves global existence for smooth initial data. The result
for general data in the Sobolev class is obtained by a standard approximation
argument.

\section{Preliminaries}

\subsection{Local well-posedness and blow-up criterion}

We recall the standard local well-posedness theorem for strictly hyperbolic
quasilinear wave equations. The following proposition can be shown by applying a general result of the local well-posedness (e.g. Hughes, Kato and Marsden \cite{HKM} and Taylor \cite{Taylor1991})

\begin{proposition}[Local well-posedness and blow-up criterion]\label{prop:local}
Let \(s>1/2\). Assume \eqref{ass:c}. For any
\[
  (u_0,u_1)\in H^{s+1}(\R)\times H^s(\R),
\]
there exist \(T>0\) and a unique solution \(u\) to
\eqref{eq:main}--\eqref{eq:initial} such that
\begin{equation}\label{eq:local-class}
  u\in \bigcap_{j=0}^2 C^j([0,T];H^{s+1-j}(\R)).
\end{equation}

Moreover, the following standard energy estimate holds. For every
\(0<T_0<T\), there exists a positive constant
\[
  C_{T_0}
  =
  C\left(
    T_0,
    \sup_{0\leq t\leq T_0}
    \bigl(
      \|u_t(t)\|_{L^\infty}
      +
      \|u_x(t)\|_{L^\infty}
    \bigr)
  \right)
\]
such that
\begin{equation}\label{eq:energy-estimate}
  \|u_t(t)\|_{H^s}
  +
  \|u_x(t)\|_{H^s}
  \leq
  C_{T_0}
  \left(
    \|u_1\|_{H^s}
    +
    \|u_0'\|_{H^s}
  \right),
  \qquad 0\leq t\leq T_0 .
\end{equation}

Let \(T^*\) be the maximal existence time of the solution in the class
\eqref{eq:local-class}. If \(T^*<\infty\), then
\begin{equation}\label{eq:blowup-criterion}
  \limsup_{t\uparrow T^*}
  \left(
    \|u_t(t)\|_{L^\infty}
    +
    \|u_x(t)\|_{L^\infty}
  \right)
  =\infty .
\end{equation}
\end{proposition}

\begin{remark}
Under the assumption \(c(\theta)\geq c_*>0\), the equation remains strictly
hyperbolic as long as \(u\) is finite. Therefore, the only obstruction to
continuing the local solution in the above Sobolev class is the blow-up of
\[
  \|u_t(t)\|_{L^\infty}+\|u_x(t)\|_{L^\infty}.
\]
Therefore we can define the maximal existence time of the solution \(u\)
of the Cauchy problem \eqref{eq:main}--\eqref{eq:initial} by 
\begin{equation}\label{def:Tstar}
  T^*
  :=
  \sup\left\{
    T\geq 0 \; ; \;
    \sup_{0\leq t\leq T}
    \left(
      \|\partial_t u(t)\|_{L^\infty}
      +
      \|\partial_x u(t)\|_{L^\infty}
    \right)
    <\infty
  \right\}.
\end{equation}
\end{remark}

\subsection{Riemann variables}

We introduce the Riemann variables
\begin{equation}\label{def:RS}
  R=u_t+c(u)u_x,\qquad
  S=u_t-c(u)u_x.
\end{equation}
Then
\begin{equation}\label{eq:utux-RS}
  u_t=\frac{R+S}{2},
  \qquad
  u_x=\frac{R-S}{2c(u)}.
\end{equation}

Let
\[
  D_-:=\partial_t-c(u)\partial_x,
  \qquad
  D_+:=\partial_t+c(u)\partial_x .
\]
A direct computation gives
\begin{align}
  D_-R
  &=
  R_t-c(u)R_x
  =
  \frac{c'(u)}{2c(u)}S(R-S),
  \label{eq:R-equation}
  \\
  D_+S
  &=
  S_t+c(u)S_x
  =
  \frac{c'(u)}{2c(u)}R(S-R).
  \label{eq:S-equation}
\end{align}
Indeed, using \(u_{tt}=c(u)^2u_{xx}\), we have
\begin{align*}
  R_t-cR_x
  &=
  \{u_{tt}+c'(u)u_tu_x+c(u)u_{tx}\}
  -c(u)\{u_{tx}+c'(u)u_x^2+c(u)u_{xx}\}
  \\
  &=
  c'(u)u_x\{u_t-c(u)u_x\}
  =
  \frac{c'(u)}{2c(u)}(R-S)S.
\end{align*}
The computation for \(S\) is the same:
\begin{align*}
  S_t+cS_x
  &=
  \{u_{tt}-c'(u)u_tu_x-c(u)u_{tx}\}
  +c(u)\{u_{tx}-c'(u)u_x^2-c(u)u_{xx}\}
  \\
  &=
  -c'(u)u_x\{u_t+c(u)u_x\}
  =
  \frac{c'(u)}{2c(u)}R(S-R).
\end{align*}

\section{Proof of the main theorem}

\subsection{Proof for smooth initial data}
We first consider the case
\[
  (u_0,u_1)\in H^3(\mathbb{R})\times H^2(\mathbb{R}).
\]
Then the local solution satisfies \eqref{eq:main} in the classical sense.

Let \(u\) be the local solution given by Proposition \ref{prop:local}, and let
\(T^*\) be its maximal existence time. We prove that \(T^*=\infty\).

Set
\begin{equation}\label{def:R0S0}
  R_0(x):=u_1(x)+c(u_0(x))u_0'(x),
  \qquad
  S_0(x):=u_1(x)-c(u_0(x))u_0'(x).
\end{equation}

\subsection{Upper estimates for \(R\) and \(S\)}

We first derive pointwise upper estimates for \(R\) and \(S\).

Let \(x_\pm (\tau)\) be the characteristics for \(R\) \(s\), defined by
\begin{equation}\label{eq:minus-char}
  \frac{d}{d\tau}x_\pm (\tau)
  =
  \pm c(u(\tau,x_\pm (\tau))),
  \qquad
  x_\pm (t)=x.
\end{equation}

From  the relations that

\[
  D_-u=S,\qquad D_+u=R,
\]
we obtain
\begin{align}
  D_-\left(\frac{R}{\sqrt{c(u)}}\right)
  &=
  -\frac{c'(u)}{2c(u)^{3/2}}S^2
  \leq 0,
  \label{eq:weighted-R}
  \\
  D_+\left(\frac{S}{\sqrt{c(u)}}\right)
  &=
  -\frac{c'(u)}{2c(u)^{3/2}}R^2
  \leq 0.
  \label{eq:weighted-S}
\end{align}
Thus we have
\[
  \frac{d}{d\tau}
  \left(
    \frac{R(\tau,x_-(\tau))}
         {\sqrt{c(u(\tau,x_-(\tau)))}}
  \right)
  \leq 0.
\]
Hence
\[
  \frac{R(t,x)}{\sqrt{c(u(t,x))}}
  \leq
  \frac{R_0(x_-(0))}{\sqrt{c(u_0(x_-(0)))}}.
\]
Using \(c_*\leq c\leq c^*\), we obtain
\begin{equation*}
  R(t,x)
  \leq
  P_R,
\end{equation*}
where
\begin{equation*}
  P_R
  =
  \sqrt{\frac{c^*}{c_*}}
  \| (R_0)_+\|_{L^\infty}.
\end{equation*}

Similarly, we have that
\begin{equation*}
  S(t,x)
  \leq
  P_S,
\end{equation*}
where
\begin{equation*}
  P_S
  =
  \sqrt{\frac{c^*}{c_*}}
  \| (S_0)_+\|_{L^\infty}.
\end{equation*}

In particular, putting
\begin{equation*}
  P=\max\{P_R,P_S\},
\end{equation*}
we have
\begin{equation}\label{eq:upper-common}
  R(t,x)\leq P,\qquad S(t,x)\leq P
\end{equation}
for all \((t,x)\in [0,T^*)\times\R\).

\subsection{Lower estimates for \(R\) and \(S\)}

Next we prove a lower estimate.  we define
\begin{equation*}
  A=\left\|\frac{c'}{2c}\right\|_{L^\infty(\R)}.
\end{equation*}
By \eqref{ass:c} and \eqref{ass:cp}, $A$ is finite.

Let
\begin{equation}\label{def:m0}
  m_0:=
  \min\left\{
    0,\inf_{x\in\R}R_0(x),\inf_{x\in\R}S_0(x)
  \right\}.
\end{equation}
We claim that
\begin{equation}\label{eq:lower-claim}
  R(t,x)\geq y(t),
  \qquad
  S(t,x)\geq y(t),
\end{equation}
where \(y(t)\) is the solution to
\begin{equation}\label{eq:y-ode}
  y'(t)=A\{Py(t)-P^2\},
  \qquad
  y(0)=m_0.
\end{equation}
Explicitly,
\begin{equation}\label{eq:y-explicit}
  y(t)=P+(m_0-P)e^{APt}.
\end{equation}

We prove \eqref{eq:lower-claim} by a comparison argument. Fix \(\eta>0\) and let \(y_\eta(t)\) be the solution of
\begin{equation}\label{eq:yeta}
  y_\eta'(t)=A\{Py_\eta(t)-P^2\}-\eta,
  \qquad
  y_\eta(0)=m_0-\eta.
\end{equation}
It is enough to show that
\begin{equation}\label{eq:yeta-lower}
  R(t,x)>y_\eta(t),
  \qquad
  S(t,x)>y_\eta(t)
\end{equation}
on every compact time interval contained in \([0,T^*)\). Then \(\eta\downarrow0\) gives \eqref{eq:lower-claim}.

At \(t=0\), by the definition of \(m_0\),
\[
  R_0(x)\geq m_0>m_0-\eta=y_\eta(0),
  \qquad
  S_0(x)\geq m_0>m_0-\eta=y_\eta(0).
\]
Suppose, for contradiction, that \eqref{eq:yeta-lower} fails. Let \(t_0 \in (0,T^*)\) be the first time such that there exists \(x_0\in\R\) with
\[
  R(t_0,x_0)=y_\eta(t_0)
  \quad \text{or} \quad
  S(t_0,x_0)=y_\eta(t_0).
\]
Since \(R(t,\cdot)\) and \(S(t,\cdot)\) are continuous and vanish at spatial
infinity, and since \(y_\eta(t)<0\), such a first point of equality is attained
at some finite point \(x_0\in\R\), if the inequality fails.

We treat the first case that \( R(t_0,x_0)=y_\eta(t_0)\). Let \(X(\tau)\) be the characteristic defined by
\[
  \frac{dX}{d\tau}=-c(u(\tau,X(\tau))),
  \qquad
  X(t_0)=x_0.
\]
Put
\[
  \Phi(\tau)=R(\tau,X(\tau))-y_\eta(\tau).
\]
By the definition of \(t_0\),
\[
  \Phi(\tau)>0\quad (\tau<t_0),
  \qquad
  \Phi(t_0)=0.
\]
Hence the derivative at \(t_0\) satisfies
\begin{equation}\label{eq:left-derivative}
  \Phi'(t_0)\leq 0.
\end{equation}
On the other hand, we have
\[
  \Phi'(t_0)
  =
  (D_-(R-y_\eta))(t_0,x_0).
\]
Since \(R(t_0,x_0)=y_\eta(t_0)\), we have
\[
  D_-(R-y_\eta)
  =
  \frac{c'(u)}{2c(u)}S(y_\eta-S)-y_\eta'
\]
at \((t_0,x_0)\).

At the  time $t_0$, we still have
\[
  S(t_0,x_0)\geq y_\eta(t_0).
\]
Moreover, by \eqref{eq:upper-common},
\[
  S(t_0,x_0)\leq P.
\]
Therefore we have
\[
  y_\eta(t_0)\leq S(t_0,x_0)\leq P.
\]
For any \(S\in[y_\eta,P]\), we have
\begin{equation}\label{eq:elementary-ineq}
  S(y_\eta-S)\geq Py_\eta-P^2.
\end{equation}
Indeed, it follows that
\[
  S(y_\eta-S)-(Py_\eta-P^2)
  =
  (P-S)(P+S-y_\eta)\geq0.
\]
Since \(Py_\eta-P^2\leq0\) and \(0\leq c'(u)/(2c(u))\leq A\), we obtain
\[
  \frac{c'(u)}{2c(u)}S(y_\eta-S)
  \geq
  A(Py_\eta-P^2).
\]
Using \eqref{eq:yeta}, we conclude that
\[
  D_-(R-y_\eta)(t_0,x_0)
  \geq
  A(Py_\eta-P^2)
  -
  (A(Py_\eta-P^2)-\eta)
  =
  \eta>0.
\]
This contradicts \eqref{eq:left-derivative}.

The case where \(S(t_0,x_0)=y_\eta(t_0)\) is treated in the same way.
Let \(X(\tau)\) be the characteristic defined by
\[
  \frac{dX}{d\tau}=c(u(\tau,X(\tau))),
  \qquad
  X(t_0)=x_0.
\]
Then
\[
  \Psi(\tau):=S(\tau,X(\tau))-y_\eta(\tau)
\]
satisfies \(\Psi(\tau)>0\) for \(\tau<t_0\), \(\Psi(t_0)=0\), and hence
\[
  \Psi'(t_0)\leq0.
\]
On the other hand, we have
\[
  D_+(S-y_\eta)
  =
  \frac{c'(u)}{2c(u)}R(y_\eta-R)-y_\eta'.
\]
Since
\[
  y_\eta(t_0)\leq R(t_0,x_0)\leq P,
\]
the same elementary inequality gives
\[
  R(y_\eta-R)\geq Py_\eta-P^2.
\]
Therefore we have
\[
  D_+(S-y_\eta)(t_0,x_0)\geq \eta>0,
\]
which is again a contradiction.

Thus \eqref{eq:yeta-lower} holds. Letting \(\eta\downarrow0\), we obtain
\eqref{eq:lower-claim}.

Consequently, for every \(T<T^*\),
\begin{equation}\label{eq:RS-bounded}
  P+(m_0-P)e^{APT}
  \leq
  R(t,x),S(t,x)
  \leq
  P
\end{equation}
for all \((t,x)\in[0,T]\times\R\).

\subsection{Completion of the proof for smooth initial data}

From \eqref{eq:utux-RS}, \eqref{eq:upper-common}, and \eqref{eq:RS-bounded}, we obtain, for every \(T<T^*\),
\begin{equation}\label{eq:utux-bound}
  \sup_{0\leq t\leq T}
  \left(
    \|u_t(t)\|_{L^\infty}
    +
    \|u_x(t)\|_{L^\infty}
  \right)
  <\infty.
\end{equation}
Indeed, since \(c(u)\geq c_*>0\),
\[
  |u_t|
  \leq
  \frac{|R|+|S|}{2},
  \qquad
  |u_x|
  \leq
  \frac{|R|+|S|}{2c_*}.
\]
Thus the \(L^\infty\)-norms of \(u_t\) and \(u_x\) are bounded on every finite time interval.

Suppose that \(T^*<\infty\). Taking \(T<T^*\) and then letting \(T\uparrow T^*\) in \eqref{eq:utux-bound}, we get
\[
  \sup_{0\leq t<T^*}
  \left(
    \|u_t(t)\|_{L^\infty}
    +
    \|u_x(t)\|_{L^\infty}
  \right)
  <\infty.
\]
This contradicts the blow-up criterion \eqref{eq:blowup-criterion}. Therefore \(T^*=\infty\), and the proof of Theorem \ref{thm:main} is complete for smooth initial data.

\subsection{Approximation argument}

We now prove the theorem for general initial data
\[
  (u_0,u_1)\in H^{s+1}(\mathbb{R})\times H^s(\mathbb{R}),
  \qquad s>\frac12 .
\]
Let \(\rho\in C_0^\infty(\mathbb{R})\) be a standard mollifier and set
\[
  \rho_\varepsilon(x)
  =
  \frac1\varepsilon \rho\left(\frac{x}{\varepsilon}\right).
\]
We define
\[
  u_0^\varepsilon:=\rho_\varepsilon*u_0,
  \qquad
  u_1^\varepsilon:=\rho_\varepsilon*u_1.
\]
Then
\[
  (u_0^\varepsilon,u_1^\varepsilon)
  \to
  (u_0,u_1)
  \quad
  \text{in }
  H^{s+1}(\mathbb{R})\times H^s(\mathbb{R})
\]
as \(\varepsilon\to0\), and
\[
  (u_0^\varepsilon,u_1^\varepsilon)\in H^3(\mathbb{R})\times H^2(\mathbb{R})
\]
for each \(\varepsilon>0\).

Let \(u^\varepsilon\) be the global solution corresponding to the initial
data \((u_0^\varepsilon,u_1^\varepsilon)\). We also set
\[
  R^\varepsilon
  =
  \partial_tu^\varepsilon+c(u^\varepsilon)\partial_xu^\varepsilon,
  \qquad
  S^\varepsilon
  =
  \partial_tu^\varepsilon-c(u^\varepsilon)\partial_xu^\varepsilon.
\]
By the estimates obtained in the previous subsection, \(R^\varepsilon\) and
\(S^\varepsilon\) are bounded from above and below uniformly in
\(\varepsilon\) on every finite time interval. 
Indeed, by Sobolev's embedding \(H^s(\R)\hookrightarrow L^\infty(\R)\) for
\(s>1/2\) and the standard properties of mollifiers, the quantities
\[
  \|(R_0^\varepsilon)_+\|_{L^\infty},\quad
  \|(S_0^\varepsilon)_+\|_{L^\infty},\quad
  -\inf_{\R}R_0^\varepsilon,\quad
  -\inf_{\R}S_0^\varepsilon
\]
are bounded uniformly in \(\varepsilon\) by a constant depending only on
\(\|u_0\|_{H^{s+1}}\) and \(\|u_1\|_{H^s}\).

More precisely, for every
\(T>0\), there exists a constant \(C_T>0\), independent of \(\varepsilon\),
such that
\[
  \sup_{0\leq t\leq T}
  \left(
    \|R^\varepsilon(t)\|_{L^\infty}
    +
    \|S^\varepsilon(t)\|_{L^\infty}
  \right)
  \leq C_T .
\]
Consequently,
\[
  \sup_{0\leq t\leq T}
  \left(
    \|\partial_tu^\varepsilon(t)\|_{L^\infty}
    +
    \|\partial_xu^\varepsilon(t)\|_{L^\infty}
  \right)
  \leq C_T .
\]

By the energy estimate in Proposition \ref{prop:local}, we also obtain the
uniform Sobolev bound
\[
  \sup_{0\leq t\leq T}
  \left(
    \|u^\varepsilon(t)\|_{H^{s+1}}
    +
    \|\partial_tu^\varepsilon(t)\|_{H^s}
  \right)
  \leq C_T
  \left(
    \|u_0\|_{H^{s+1}}
    +
    \|u_1\|_{H^s}
  \right),
\]
where \(C_T\) is independent of \(\varepsilon\).

Let \(u\) be the local solution corresponding to the initial data
\((u_0,u_1)\), and let \(T^*\) be its maximal existence time. By the continuous
dependence of the local solution on the initial data, for every \(T<T^*\),
\[
  \sup_{0\leq t\leq T}
  \left(
    \|u^\varepsilon(t)-u(t)\|_{H^1}
    +
    \|\partial_tu^\varepsilon(t)-u_t(t)\|_{L^2}
  \right)
  \to0
\]
as \(\varepsilon\to0\).

Moreover, by the lower semicontinuity of the Sobolev norm, we obtain
\[
  \sup_{0\leq t\leq T}
  \left(
    \|u(t)\|_{H^{s+1}}
    +
    \|u_t(t)\|_{H^s}
  \right)
  \leq
  C_T
  \left(
    \|u_0\|_{H^{s+1}}
    +
    \|u_1\|_{H^s}
  \right)
\]
for every \(T<T^*\). In particular, the solution \(u\) cannot blow up at
finite time in the sense of the blow-up criterion in Proposition
\ref{prop:local}. Therefore \(T^*=\infty\).

This proves the global existence part of Theorem \ref{thm:main}. The
uniqueness and continuous dependence follow from the local well-posedness
theory and the same a priori estimates. Hence the proof of Theorem
\ref{thm:main} is complete.

\bibliographystyle{amsplain}
\bibliography{bibliography}
\end{document}